\documentclass{article}
\usepackage{amssymb}
\usepackage{graphicx}
\usepackage{amsmath}

\newtheorem{theorem}{Theorem}

\newtheorem{corollary}[theorem]{Corollary}

\newtheorem{definition}[theorem]{Definition}

\newtheorem{lemma}[theorem]{Lemma}

\newtheorem{proposition}[theorem]{Proposition}

\input{tcilatex}
\begin{document}

\title{Klein geometries, parabolic geometries \\
and differential equations of finite type}
\author{Ender Abado\u{g}lu, Erc\"{u}ment Orta\c{c}gil, Ferit \"{O}zt\"{u}rk}
\maketitle

\begin{abstract}
We define the infinitesimal and geometric orders of an effective Klein
geometry $G/H.$ Using these concepts, we prove $i)$ For any integer $m\geq
2, $ there exists an effective Klein geometry $G/H$ of infinitesimal order $%
m $ such that $G/H$ is a projective variety (Corollary 9). $ii)$ An
effective Klein geometry $G/H$ of geometric order $M$ defines a differential
equation of order $M+1$ on $G/H$ whose global solution space is $G$
(Proposition 18).
\end{abstract}

\section{Introduction}

As stated in the review of the book [19], in current mathematics it is
mainly the scores of Lie groups and Lie algebras that are extracted from the
great symphony of Sophus Lie. Later, mainly due to the efforts of E.Cartan,
W.Killing and H.Weyl, the classification of semisimple Lie algebras and
their representations is achieved as one the most beautiful and complete
theories in mathematics. Going back to the transformation groups of S.Lie,
what do the root systems, Dynkin diagrams, Weyl groups etc. really
correspond to in terms of the infinitesimal generators of a transitive and
effective action of the Lie group? Conversely, if we start with a transitive
and effective action of a Lie group, what do prolongations, differential
invariants and other fundamental concepts arising in Lie's framework of
transformation groups correspond to in terms of the Lie group and its Lie
algebra? The purpose of this paper is to take a modest step towards the
answers of these questions in the spirit of the framework proposed in [13].
The technical results stated in the above abstract appear as byproducts.

This paper is organized as follows. In Section 2 we work in the ``universal
envelope'' of [13] and recall some well known facts about the infinite
dimensional Lie algebra of formal vector fields $J_{\infty }\mathfrak{X}_{p}$
with bracket $\{$ $,$ $\}_{\infty }$ and $k$-jets of vector fields $J_{k}%
\mathfrak{X}_{p}$ with the algebraic bracket $\{$ $,$ $\}_{k}$ induced by $\{
$ $, $ $\}_{\infty },$ refering to [4], [8], [17], [14] for more details. $\{
$ $,$ $\}_{k}$ reduces the order by one (see (3)) and therefore does \textit{%
not }endow $J_{k}\mathfrak{X}_{p}$ with a Lie algebra structure.

In Section 3, we restrict our attention to geometries contained in the
universal filtration (1) and are defined by effective Klein geometries $G/H.$
These geometries define filtrations which stabilize at zero after a finite
number of steps (see (9)). Using this fact, we define the infinitesimal
order of an effective Klein geometry. This concept exists in the fundamental
works [8], [17] and also in [4], pg.6. Now $\{$ $,$ $\}_{k}$ restricts to
the $k$-jets of the infinitesimal generators $J_{k}(\overline{\mathfrak{g}}%
)_{o}\subset J_{k}\mathfrak{X}_{p}$ of the action of $G$ on $G/H$ at some $%
o\in G/H.$ On the other hand, we define another bracket $\ [$ $,$ $]_{k}$
(see (11)) which also reduces the order by one. The definition of this new
bracket uses only the Lie algebras $\mathfrak{g}$, $\mathfrak{h}$ of $G$ and 
$H$ and seems to be unrelated to jets. Our main result (Proposition 5) shows
that $\{ $ $,$ $\}_{k}$ and $[$ $,$ $]_{k}$ coincide and become the bracket
of $\mathfrak{g}$ if $G$ acts effectively on $G/H$ and $k$ is sufficiently
large. This result allows us to detect jets inside $(G,H)$ using only group
theory.

In Section 4, we consider $\left| k\right| $-graded semisimple Lie algebras.
This well known grading is induced by a choice of positive simple roots and
seems to be totally unrelated to the grading in terms of jets used before.
However, the main result of Section 4 (Proposition 8) shows that they indeed
coincide. This fact implies the first result in the above abstract
(Corollary 9), gives an affirmative answer to a question in [13] and also
settles an open problem posed in [20] on pg. 325. It also opens the way to
express many standard concepts in the theory of semisimple Lie groups in
terms of jet theory by expressing them in terms of the coefficients of the
Taylor expansions of the infinitesimal generators, but much remains to be
done in this direction.

In Section 5, we introduce the concept of geometric order of an effective
Klein geometry. This concept exists also in [20] and is implicit in [16]. We
prove the anolog of Proposition 5 on the group level (Proposition 14) and
derive some consequences. The most notable is Corollary 16 which shows the
existence of some canonical splittings.

In Section 6 we show that Corollary 17 implies the second result in the
abstract. The mentioned differential equation is an $ODE$ if $\dim G/H=1$
and a system of $PDE$'s if $\dim G/H\geq 2.$ This differential equation
reduces to Lie's First fundamental Theorem when geometric order is zero and
also generalizes the well known Schwarzian differential equation for Mobius
transformations.

In the Appendix we make some comments on the relation of our work to [8],
[17], [6].

\section{Formal vector fields}

In this section we recall some well known facts in the form which will be
needed in the next sections. We refer to [4], [8], [17], [14] for more
details.

Let $M$ be a differentiable manifold with $\dim M=n$ and $p\in M.$ Let $%
\mathfrak{X}$ denote the Lie algebra of smooth vector fields on $M$ and for $%
X\in \mathfrak{X}$, let $(j_{\infty }X)_{p}$ denote the $\infty $-jet of $X$
at $p.$ We define the vector space $(J_{\infty }\mathfrak{X})_{p}\doteq
\{(j_{\infty }X)_{p}\mid X\in \mathfrak{X}\}.$ For simplicity of notation,
we denote $(J_{\infty }\mathfrak{X})_{p}$ by $J_{\infty }$ in this section.
Let $\widetilde{J}_{k}\subset J_{\infty }$ denote the subspace consisting of
those $(j_{\infty }X)_{p}$ vanishing at all orders up to and including order 
$k\geq 0.$ We set $\widetilde{J}_{-1}\doteq J_{\infty }.$ Thus we obtain the
following descending filtration of subspaces

\begin{equation}
...\subset \widetilde{J}_{2}\subset \widetilde{J}_{1}\subset \widetilde{J}%
_{0}\subset \widetilde{J}_{-1}=J_{\infty }
\end{equation}

In the spirit of the framework proposed in [13], we call (1) the universal
filtration at $p\in M.$

We now define the vector space $\widehat{J}_{k}\doteq \widetilde{J}_{k}/$ $%
\widetilde{J}_{k+1}$ $,-1\leq k$. Note that $\widehat{J}_{-1}=T(M)_{p}=$ the
tangent space of $M$ at $p.$ Thus we obtain

\begin{equation}
J_{\infty }=\widehat{J}_{-1}\oplus \widehat{J}_{0}\oplus \widehat{J}%
_{1}\oplus \widehat{J}_{2}\oplus .....
\end{equation}
We define a bracket $\{$ $,$ $\}_{\infty }$ on $J_{\infty }$ by $%
\{(j_{\infty }X)_{p},(j_{\infty }Y)_{p}\}_{\infty }\doteq (j_{\infty
}[X,Y])_{p}.$ This gives a Lie algebra homomorphism $\mathfrak{X}\rightarrow
J_{\infty }$ defined by $X\rightarrow (j_{\infty }X)_{p}.$ The bracket $\{$ $%
,$ $\}_{\infty }$ turns $\widehat{J}_{0}\oplus \widehat{J}_{1}\oplus 
\widehat{J}_{2}\oplus ...$ into a graded Lie algebra: $\{\widehat{J}_{i},%
\widehat{J}_{j}\}\subset \widehat{J}_{i+j}$, $0\leq i,j.$ We also have $[%
\widehat{J}_{-1},\widehat{J}_{i}]\subset \widehat{J}_{i-1},$ $i\geq 0$ but $%
[x,y]$ is undefined for $x,y\in \widehat{J}_{-1}$ which can be checked using
coordinates. It is standard to define $[\widehat{J}_{-1},\widehat{J}_{-1}]=0$
and turn $J_{\infty }$ into a graded Lie algebra by setting $\widehat{J}%
_{-2}\doteq 0.$ This definition turns out to be incompatible with the
present framework (see the paragraph below (24)). For this reason, we leave $%
[\widehat{J}_{-1},\widehat{J}_{-1}]$ undefined. Note that (1) is now a
descending filtration of ideals inside \textit{\ }$J_{0}$ but not inside%
\textit{\ }$\widetilde{J}_{-1}=J_{\infty }$ since $[\widetilde{J}_{-1},%
\widetilde{J}_{k}]\subset \widetilde{J}_{k-1}.$

We now truncate (2) at $k-1$ and define $J_{k}\doteq $ $\widehat{J}%
_{-1}\oplus \widehat{J}_{0}\oplus .....\oplus \widehat{J}_{k-1},$ $0\leq k,$
so that $J_{0}=\widehat{J}_{-1}$. Clearly, $J_{k}=J_{\infty }/\widetilde{J}%
_{k}.$ An element of $J_{k}$\ is called a $k$-jet of a vector field at $p$
and is denoted by $(j_{k}X)_{p}.$ Thus $(j_{k}X)_{p}=[(j_{\infty
}X)_{p}]_{k} $ where $[(j_{\infty }X)_{p}]_{k}$ denotes the equivalence
class of $(j_{\infty }X)_{p}$ in $J_{\infty }/\widetilde{J}_{k}$.

The bracket $\{$ $,$ $\}_{\infty }$ gives the algebraic bracket

\begin{equation}
\{\text{ },\text{ }\}_{k}:\quad J_{k}\times J_{k}\rightarrow J_{k-1}\qquad
1\leq k
\end{equation}
defined as follows: For $(j_{k}X)_{p}$, $(j_{k}Y)_{p}$ $\in J_{k}$ where $%
(j_{k}X)_{p}=[(j_{\infty }X)_{p}]_{k}$ and $(j_{k}Y)_{p}=[(j_{\infty
}Y)_{p}]_{k},$ we define $\{(j_{k}X)_{p},(j_{k}Y)_{p}\}_{k}\doteq \lbrack
(j_{\infty }[X,Y])_{p}]_{k-1}.$

Let $J_{k,j}$ denote the kernel of the projection map $\pi
_{k,j}:J_{k}\rightarrow J_{j},$ $0\leq $ $j\leq k-1.$ Thus we have the exact
sequence

\begin{equation}
0\longrightarrow J_{k,j}\longrightarrow J_{k}\overset{\pi _{k,j}}{%
\longrightarrow }J_{j}\longrightarrow 0
\end{equation}

Now $\{$ $,$ $\}_{k}$ restricts to $J_{k,0}$ as

\begin{equation}
\{\text{ },\text{ }\}_{k}:\quad J_{k,0}\times J_{k,0}\rightarrow J_{k,0}
\end{equation}

Thus $J_{k,0}$ is a Lie algebra with bracket $\{$ $,$ $\}_{k}$ and is called
the isotropy subalgebra. In fact, let $(\mathfrak{G}_{k})_{p}^{p}$ be the
Lie group of $k$-jets of local diffeomorphisms with source and target at $p.$
Any choice of coordinates near $p$ identifies $(\mathfrak{G}_{k})_{p}^{p}$
with the $k$'th order jet group $GL_{k}(n)$ and $J_{k,0}$ is the Lie algebra
of $(\mathfrak{G}_{k})_{p}^{p}.$

Now let $(\mathfrak{G}_{k})_{q}^{p}$ denote the set of all $k$-jets of local
diffeomorphisms with source at $p\in M$ and target at $q\in M.$ If $%
j_{k+1}(f)_{q}^{p}\in (\mathfrak{G}_{k+1})_{q}^{p}$, then $j_{k+1}(f)_{q}^{p}
$ induces an isomorphism $\natural j_{k+1}(f)_{q}^{p}:(J_{k}\mathfrak{X}%
)_{p}\rightarrow (J_{k}\mathfrak{X})_{q}$. In particular, we obtain a
representation of $(\mathfrak{G}_{k+1})_{p}^{p}$ on $J_{k},$ that is, a
homomorphism

\begin{equation}
\natural :(\mathfrak{G}_{k+1})_{p}^{p}\longrightarrow GL(J_{k})
\end{equation}
defined by $j_{k+1}(f)_{p}^{p}\rightarrow \natural j_{k+1}(f)_{p}^{p}.$ The
representation $\natural $ \ is faithful.

For an explicit formula for $\natural $ to be used in Section 5$,$ let $%
(j_{k}X)_{p}\in J_{k}$ and $j_{k+1}(f)_{p}^{p}\in (\mathfrak{G}%
_{k+1})_{p}^{p}.$ The diffeomorphism $f\circ e^{tX}\circ f^{-1}$ is the
identity when $t=0$ and defines curves starting at points near $p.$
Differentiating these curves at $t=0,$ we obtain a vector field defined near 
$p$ which we denote by $\frac{d}{dt}(f\circ e^{tX}\circ f^{-1})_{\mid t=0}.$
We have

\begin{equation}
\natural j_{k+1}(f)_{p}^{p}((j_{k}X)_{p})=j_{k}(\frac{d}{dt}(f\circ
e^{tX}\circ f^{-1})_{\mid t=0})_{p}
\end{equation}

The infinitesimal representation induced by (6) is

\begin{equation}
\flat :J_{k+1,0}\longrightarrow gl(J_{k})
\end{equation}

A subset $\mathfrak{T}\subset J_{k}$ is called transitive if $\pi _{k,0}(%
\mathfrak{T})=J_{k,0}$ $\ =T(M)_{p}.$

The proof the next lemma is a straightforward computation in local
coordinates.

\begin{lemma}
Let $\mathfrak{T}\subset J_{k}$ be transitive and $j_{k+1}(f)_{p}^{p}\in (%
\mathfrak{G}_{k+1})_{p}^{p}$. If $\natural j_{k+1}(f)_{p}^{p}(Y)=Y$ for all $%
Y\in $ $\mathfrak{T}$, then $j_{k+1}(f)_{p}^{p}=id$.
\end{lemma}

The local formulas for $\{$ $,$ $\}_{\infty }$ are obtained by
differentiating the usual bracket formula $[X,Y]=X^{a}\partial
_{a}Y^{i}-Y^{a}\partial _{a}X^{i}$ successively infinitely many times,
evaluating at $p$ and substituting jets. We denote the formula obtained by $%
k $-times differentiation by $A_{k},$ $k\geq 0.$ For instance, $A_{0}$ is $%
\{(j_{\infty }X)_{p},(j_{\infty }Y)_{p}\}^{i}=X^{a}Y_{a}^{i}-Y^{a}X_{a}^{i}.$
We make use of these formulas in the proof of Lemma 4.

It is also crucial to observe that if we replace $\mathfrak{X}$ in the above
construction by the germs of vector fields at $p,$ we get the same $%
J_{\infty }$, since a partition of unity argument shows that any such germ
comes from some $X\in \mathfrak{X}$. This is far from being true in the
(complex) analytic category. In this case, if $(j_{\infty }X)_{p}=(j_{\infty
}Y)_{p}$ for some $X,Y\in \mathfrak{X}$, then $X=Y$ on $M$ and so $J_{\infty
}$ contains global information. We have a special case of this situation
below where for some integer $m,$ $(j_{m}X)_{p}$ uniquely determines $X$ and
(3) turns into an honest Lie bracket on some finite dimensional Lie algebra
for $k=m+1.$

\section{Infinitesimal order}

Let $G$ be a Lie group (not necessarily connected) and $H$ a Lie subgroup. $%
G $ acts on $G/H$ by $L_{g}(xH)=gxH.$ Other than this action, there is
another fundamental concept inherent in the definition of the homogeneous
space $G/H: $ the $H$-principal bundle $G\rightarrow G/H.$ To emphasize our
choice in this paper, we make the following

\begin{definition}
A global Klein geometry consists of the following:

$i)$ A homogeneous space $G/H$

$ii)$ The (left) action of $G$ on $G/H$
\end{definition}

In this paper, Klein geometry means a global Klein geometry. We denote a
Klein geometry by $G/H$ and call $G/H$ effective if $G$ acts effectively. If 
$K$ is the largest normal subgroup of $G$ contained in $H$, then $G/H$ is
effective iff $K=\{e\}.$

We denote $H$ by $G_{0}$ and the Lie algebras of $G,$ $G_{0}$ by $\mathfrak{g%
},$ $\mathfrak{g}_{0}.$ Now following [8], [17], we inductively define $%
\mathfrak{g}_{k+1}\doteq \{x\in \mathfrak{g}_{k}\mid \lbrack x,\mathfrak{g}%
]\subset \mathfrak{g}_{k}, $ $0\leq k\}.$ Then $\mathfrak{g}_{k+1}\subset 
\mathfrak{g}_{k}$ is an ideal for $0\leq k.$ Since $\mathfrak{g}$ is finite
dimensional, there exists an integer $m$ such that $\mathfrak{g}_{m+i}=%
\mathfrak{g}_{m}$ for all $0\leq i$. Since $\mathfrak{g}_{m}\subset 
\mathfrak{g}$ is an ideal, there exists a connected and normal subgroup $%
H\vartriangleleft G_{-1}$ with Lie algebra $\mathfrak{g}_{m}.$ Now $H$ is
contained in the connected component of $G_{0}$ since $\mathfrak{g}%
_{m}\subset \mathfrak{g}_{0}$. Therefore $H=\{e\}$ if $G$ acts effectively
on $G/G_{0}$. So we conclude $\mathfrak{g}_{m}=0$ in this case.

\begin{definition}
The integer $m$ such that $\mathfrak{g}_{m}=0$ but $\mathfrak{g}_{m-1}\neq 0$
is called the infinitesimal order of the effective Klein geometry $G/G_{0}$.
\end{definition}

Until Section 5, order means infinitesimal order.

Therefore, an effective Klein geometry determines the descending filtration 
\begin{equation}
\{0\}\subset \mathfrak{g}_{m-1}\subset .....\subset \mathfrak{g}_{2}\subset 
\mathfrak{g}_{1}\subset \mathfrak{g}_{0}\subset \mathfrak{g}
\end{equation}

Now, since $\mathfrak{g}_{k}\subset \mathfrak{g}_{0}$ is an ideal, we have
the homomorphism $ad_{k}:\mathfrak{g}_{0}\rightarrow gl(\mathfrak{g}/%
\mathfrak{g}_{k})$ defined by $ad_{k}(x)(y+\mathfrak{g}_{i})=[x,y]+\mathfrak{%
g}_{k}=ad_{x}(y)+\mathfrak{g}_{k}$ where $[$ $,$ $]$ is the bracket of $%
\mathfrak{g}$. We observe that $\ker (ad_{k})=\mathfrak{g}_{k+1}.$ This
gives an alternative and more conceptual definition of the spaces in (9). In
particular, we obtain the faithful representation

\begin{equation}
ad_{k}:\mathfrak{g}_{0}/\mathfrak{g}_{k+1}\longrightarrow gl(\mathfrak{g}/%
\mathfrak{g}_{k})
\end{equation}
(for simplicity of notation, we keep the same notation for maps when we pass
to quotients or make identifications).

We define the bracket

\begin{equation}
\lbrack \text{ },\text{ }]_{k+1}:\quad \mathfrak{g}/\mathfrak{g}_{k+1}\times 
\mathfrak{g}/\mathfrak{g}_{k+1}\rightarrow \mathfrak{g}/\mathfrak{g}%
_{k}\qquad 0\leq k
\end{equation}
by $[a+\mathfrak{g}_{k+1},b+\mathfrak{g}_{k+1}]_{k+1}\doteq \lbrack a,b]+%
\mathfrak{g}_{k}.$ Since $[\mathfrak{g},\mathfrak{g}_{k+1}]\subset \mathfrak{%
g}_{k},$ $[$ $,$ $]_{k+1}$ is well defined. Note that $[$ $,$ $]_{m+1}$ $=[$ 
$,$ $].$ We also have the projection map

\begin{equation}
\overline{\pi }_{k,j}:\mathfrak{g}/\mathfrak{g}_{k}\rightarrow \mathfrak{g}/%
\mathfrak{g}_{j}\qquad j+1\leq k
\end{equation}
with kernel $\mathfrak{g}_{j}/\mathfrak{g}_{k}$ and the restricted bracket

\begin{equation}
\lbrack \text{ },\text{ }]_{k}:\quad \mathfrak{g}_{0}/\mathfrak{g}_{k}\times 
\mathfrak{g}_{0}/\mathfrak{g}_{k}\rightarrow \mathfrak{g}_{0}/\mathfrak{g}%
_{k}\qquad 1\leq k
\end{equation}
which is well defined since $\mathfrak{g}_{k}\subset \mathfrak{g}_{0}$ is an
ideal.

Now we also have the Lie algebra homomorphism $\mathfrak{g}\rightarrow 
\mathfrak{X}$ which maps $X\in \mathfrak{g}$ to its infinitesimal generator $%
\overline{X}\in \mathfrak{X}$, where $\mathfrak{X}$ is the Lie algebra of
smooth vector fields on $M=G/G_{0}$. We denote the Lie subalgebra of
infinitesimal generators by $\overline{\mathfrak{g}}\subset \mathfrak{X}$.
Now $\overline{\mathfrak{g}}$ is isomorphic to $\mathfrak{g}$ since the
action is both transitive and effective. We define $J_{k}(\overline{%
\mathfrak{g}})_{o}\doteq \{(j_{k}\overline{X})_{o}\mid \overline{X}\in 
\overline{\mathfrak{g}},$ $0\leq \,k\}\subset $ $(J_{k}\mathfrak{X})_{o}.$
All the constructions of Section 2 can be done now with $\overline{\mathfrak{%
g}}$ at $o$ and all the spaces obtained in this way imbed in the spaces in
Section 2 together with their grading. In particular, we obtain a filtration
contained in the universal filtration (1).

The restriction of (3) gives

\begin{equation}
\{\text{ },\text{ }\}_{k+1}:\quad J_{k+1}(\overline{\mathfrak{g}})_{o}\times
J_{k+1}(\overline{\mathfrak{g}})_{o}\rightarrow J_{k}(\overline{\mathfrak{g}}%
)_{o}\qquad 0\leq k
\end{equation}
(14) follows from [14] but can be checked also directly. Note that $J_{0}(%
\overline{\mathfrak{g}})_{o}=$ $\mathfrak{g}/\mathfrak{g}_{0}=$ $%
T(G/G_{0})_{o}.$ The restriction of $\pi _{k,j}$ gives

\begin{equation}
\pi _{k,j}:J_{k}(\overline{\mathfrak{g}})_{o}\rightarrow J_{j}(\overline{%
\mathfrak{g}})_{o}\qquad j+1\leq k
\end{equation}

Clearly (15) commutes with (14). We have

\begin{equation}
\{\text{ },\text{ }\}_{k}:\quad J_{k,0}(\overline{\mathfrak{g}})_{o}\times
J_{k,0}(\overline{\mathfrak{g}})_{o}\rightarrow J_{k,0}(\overline{\mathfrak{g%
}})_{o}\qquad 1\leq k
\end{equation}

We also have the faithful representation

\begin{equation}
\flat _{\mid (\mathfrak{g}_{-1},\mathfrak{g}_{0})}:J_{k+1,0}(\overline{%
\mathfrak{g}})_{o}\longrightarrow gl(J_{k}(\overline{\mathfrak{g}})_{o})
\end{equation}

The faithfulness of (17) follows from the infinitesimal analog of Lemma 1
and the fact that $J_{k}(\overline{\mathfrak{g}})_{o}\subset (J_{k}\mathfrak{%
X})_{o}$ is transitive.

To clarify the analogy between (11)-(14), (12)-(15), (13)-(16) and
(10)-(17), we define the map

\begin{eqnarray}
\theta _{k} &:&\mathfrak{g}\rightarrow J_{k}(\overline{\mathfrak{g}})_{o} \\
\theta _{k}(x) &=&j_{k}(\overline{x})_{o}\quad k\geq 0  \notag
\end{eqnarray}
where $\overline{x}$ is the infinitesimal generator of $x.$ Now $\theta _{k}$
is clearly linear and is surjective by the definition of $J_{k}(\overline{%
\mathfrak{g}})_{o}.$

\begin{lemma}
The kernel of $\theta _{k}$ is $\mathfrak{g}_{k}.$
\end{lemma}

Proof: We prove by induction on $k$ that, for $0\leq k$, $(j_{k}\overline{X}%
)_{o}=0$ if and only if $X\in \mathfrak{g}_{k}$.

For $k=0$, the claim is that $\overline{X}(o)=0$ if and only if $X\in 
\mathfrak{g}_{0}$. Let $\pi :G\rightarrow G/G_{0}$ be the quotient map,
inducing $\overline{\pi }=d\pi (e):\mathfrak{g}\rightarrow \mathfrak{g}/%
\mathfrak{g}_{0}$. By definition, $\overline{X}(o)=\frac{d}{dt}_{\mid
t=0}(e^{tX}(o))=\overline{\pi }(\frac{d}{dt}_{\mid t=0}e^{tX})=\overline{\pi 
}(X)$. Hence the claim follows. Now $X\in \mathfrak{g}_{k+1}$

\begin{eqnarray}
&\Longleftrightarrow &X\in \mathfrak{g}_{k}\text{ and }[X,Y]\in \mathfrak{g}%
_{k}\text{ for all }Y\in \mathfrak{g}\quad \text{(definition of \ }\mathfrak{%
g}_{k+1})  \notag \\
&\Longleftrightarrow &(j_{k}\overline{X})_{o}=0\text{ and }j_{k}\overline{%
[X,Y]}_{o}=0\text{ for all }Y\in \mathfrak{g}\quad \text{(induction
hypothesis)}  \notag \\
&\Longleftrightarrow &(j_{k}\overline{X})_{o}=0\text{ and }j_{k}[\overline{X}%
,\overline{Y}]_{o}=0\text{ for all }Y\in \mathfrak{g} \\
&\Longleftrightarrow &(j_{k}\overline{X})_{o}=0\text{ and }\{(j_{k+1}%
\overline{X})_{o},(j_{k+1}\overline{Y})_{o}\}_{k+1}=0\text{ for all }%
\overline{Y}\in J_{0}(\overline{\mathfrak{g}})_{o}\text{ }  \notag \\
&&\qquad \qquad \qquad \;\;\qquad \qquad \qquad \qquad \qquad \qquad \qquad 
\text{(definition of }\{\text{ },\text{ }\}_{k+1})\text{ }  \notag
\end{eqnarray}

We now choose some coordinate system $(x^{i})$ around $o$ and recall that
the components of $\{(j_{k+1}\overline{X})_{o},(j_{k+1}\overline{Y}%
)_{o}\}_{k+1}$ are given by the formulas $A_{0},A_{1},...,A_{k}.$ Since $%
(j_{k}\overline{X})_{o}=0$, all terms in the formulas $A_{0},A_{1},...,A_{k}$
vanish except the last term in $A_{k}$ which is $-\overline{Y}^{a}\overline{X%
}_{aj_{k}...j_{1}}^{i}.$ Therefore, the last formula in (19) is equivalent to

\begin{equation}
\overline{Y}^{a}\overline{X}_{aj_{k}...j_{1}}^{i}=0\quad \text{for all }%
\overline{Y}\in J_{0}(\overline{\mathfrak{g}})_{o}
\end{equation}

Since $J_{k}(\overline{\mathfrak{g}})_{o}$ is transitive, (20) is equivalent
to $\overline{X}_{j_{k+1}j_{k}...j_{1}}^{i}=0,$ that is, $(j_{k+1}\overline{X%
})_{o}=0.$ This completes the inductive step. \ $\square $

Lemma 4 gives the linear isomorphism

\begin{equation}
\theta _{k}:\mathfrak{g}/\mathfrak{g}_{k}\longrightarrow J_{k}(\overline{%
\mathfrak{g}})_{o}
\end{equation}
and its proof shows that (9) coincides with the filtration described above
inside the universal filtration (1).

\begin{proposition}
The diagrams
\end{proposition}

\begin{equation}
\begin{array}{ccc}
\mathfrak{g}/\mathfrak{g}_{k} & \overset{\theta _{k}}{\longrightarrow } & 
J_{k}(\overline{\mathfrak{g}})_{o} \\ 
\downarrow _{\overline{\pi }_{k,j}} &  & \downarrow _{_{\pi _{k,j}}} \\ 
\mathfrak{g}/\mathfrak{g}_{j} & \overset{\theta _{j}}{\longrightarrow } & 
J_{j}(\overline{\mathfrak{g}})_{o}%
\end{array}%
\end{equation}

\begin{equation}
\begin{array}{ccc}
\mathfrak{g}/\mathfrak{g}_{k+1}\times \mathfrak{g}/\mathfrak{g}_{k+1} & 
\overset{[\text{ },\text{ }]_{k+1}}{\longrightarrow } & \mathfrak{g}/%
\mathfrak{g}_{k} \\ 
\downarrow _{\theta _{k}\times ad_{k}} &  & \downarrow _{\theta _{k}} \\ 
J_{k+1}(\overline{\mathfrak{g}})_{o}\times J_{k+1}(\overline{\mathfrak{g}}%
)_{o} & \overset{\{\text{ },\text{ }\}_{k+1}}{\longrightarrow } & J_{k}(%
\overline{\mathfrak{g}})_{o}%
\end{array}%
\end{equation}
commute where $k\geq 0$.

Proof: The commutativity of the first diagram is straightforward. As for the
second, for $X,Y\in \mathfrak{g}_{-1}$ we have 
\begin{equation*}
\begin{array}{rl}
\{\theta _{k+1}(X+\mathfrak{g}_{k+1}),\theta _{k+1}(Y+\mathfrak{g}%
_{k+1})\}_{k+1} & =\{(j_{k+1}\overline{X})_{o},(j_{k+1}\overline{Y}%
)_{o}\}_{k+1} \\ 
& =[(j_{\infty }\{\overline{X},\overline{Y}\})_{o}]_{k} \\ 
& =[(j_{\infty }\overline{[X,Y]})_{o}]_{k} \\ 
& =(j_{k}\overline{[X,Y]})_{o}%
\end{array}%
\end{equation*}
while 
\begin{equation*}
\theta _{k}([X+\mathfrak{g}_{k+1},Y+\mathfrak{g}_{k+1}]_{k+1})=\theta
_{k}([X,Y]+\mathfrak{g}_{k})=(j_{k}\overline{[X,Y]})_{o}\qquad \square
\end{equation*}

Thus we also obtain the linear isomorphisms $\theta _{m}:\mathfrak{g}%
\rightarrow J_{m}(\overline{\mathfrak{g}})_{o},$ $\theta _{k}:\mathfrak{g}%
_{0}/\mathfrak{g}_{k}\rightarrow J_{k,0}(\overline{\mathfrak{g}})_{o}\subset
(J_{k,0}\mathfrak{X})_{o} $ and the commutative diagram

\begin{equation}
\begin{array}{cccc}
ad_{k}: & \mathfrak{g}_{0}/\mathfrak{g}_{k+1} & \longrightarrow & gl(%
\mathfrak{g}/\mathfrak{g}_{k}) \\ 
& \parallel \theta _{k} &  & \parallel \theta _{k} \\ 
\flat _{\mid (\mathfrak{g},\mathfrak{g}_{0})}: & J_{k+1,0}(\overline{%
\mathfrak{g}})_{o} & \longrightarrow & gl(J_{k}(\overline{\mathfrak{g}})_{o})%
\end{array}%
\end{equation}

Lemma 4 and Proposition 5 show that $J_{m}(\overline{\mathfrak{g}}%
)_{o}=J_{m+1}(\overline{\mathfrak{g}})_{o}$ and (11) gives the bracket $\{$ $%
,$ $\}_{m+1}:$ $J_{m}(\overline{\mathfrak{g}})_{o}\times J_{m}(\overline{%
\mathfrak{g}})_{o}\rightarrow J_{m}(\overline{\mathfrak{g}})_{o}$ which
coincides with the bracket of $\mathfrak{g}$. Therefore $\theta _{m}$ is
also a Lie algebra isomorphism. This statement holds also for $m=0$, that
is, when $G_{0}$ is discrete. For $m=0,$ note that the standard notation $%
\mathfrak{g}_{-1}$ for $\mathfrak{g}$ and the extension of the formula (11)
to $k=-1$ is incorrect in the present framework if $\mathfrak{g}$ is not
abelian.

\begin{definition}
Let $G/G_{0}$ be a any Klein geometry whose descending filtration stabilizes
at $\{0\}$. Then $G/G_{0}$ is called almost effective.
\end{definition}

If $G/G_{0}$ is almost effective, then $K$ is discrete as we show in Section
5 (see the paragraph before Lemma 10). Therefore the effective Klein
geometry $\frac{G/K}{G_{0}/K}$ defines the same filtration as $G/G_{0}$ and
has order $m.$ This fact is used in Proposition 8 and Corollary 9.

In this section we worked excusively at the point $o\in G/G_{o}.$ Using
homogeneity, all the above constructions can be done at any point $p\in
G/G_{o}$. This fact is needed only in Corollary 17.

\section{Parabolic geometries}

Let $\mathfrak{s}$ be a $\left| k\right| $-graded semisimple Lie algebra
over $\mathbb{F=R}$ or $\mathbb{C},$ $k\in \mathbb{Z}$. Thus 
\begin{equation}
\mathfrak{s}=\mathfrak{s}_{-k}\oplus \mathfrak{s}_{-k+1}\oplus ....\mathfrak{%
s}_{-1}\oplus \mathfrak{s}_{0}\oplus \mathfrak{s}_{1}\oplus ....\oplus 
\mathfrak{s}_{k-1}\oplus \mathfrak{s}_{k}
\end{equation}

To recall the grading in (25), we first assume that $\mathbb{F=C}$, $%
\mathfrak{s} $ is a semisimple Lie algebra and $\mathfrak{p}\subset 
\mathfrak{s}$ is a parabolic subalgebra. We can find a Cartan subalgebra and
a set of positive roots such that $\mathfrak{p}$ is standard with respect to
a set $\Sigma $ of simple roots. Now the $\Sigma $-height gives a $\left|
k\right| $-grading on $\mathfrak{s}$, where $k$ is the $\Sigma $-height of
the maximal root of $\mathfrak{s}$ which gives the grading (25) with $%
\mathfrak{p}=\mathfrak{s}_{0}\oplus \mathfrak{s}_{1}\oplus ....\oplus 
\mathfrak{s}_{k-1}\oplus \mathfrak{s}_{k}.$ For $\mathbb{F=R}$, the
complexification $\mathfrak{s}^{\mathbb{C}}=\mathfrak{g}$ is also $\left|
k\right| $-graded and $\mathfrak{s}$ is a real form of the complex pair $(%
\mathfrak{g},\mathfrak{p}).$ For our purpose in this section, the relevant
fact is that the origin of grading in (25) seems to be totally unrelated to
the gradings defined above in terms of jets.

Let $S$ be a Lie group with Lie algebra $\mathfrak{s}$ and $P\subset S$ a
Lie subgroup with Lie algebra $\mathfrak{p}.$ To be consistent with our
notation above, we should denote $P$ by $S_{0}$ but we will not do this.

If $x\in \mathfrak{p}$ satisfies $[x,\mathfrak{s}]\subset \mathfrak{p}$,
then the grading in (25) implies $x\in \mathfrak{s}_{k}$ and clearly $[%
\mathfrak{s}_{k},\mathfrak{s}]\subset \mathfrak{p.}$ Further, if $x\in 
\mathfrak{s}_{k}$ satisfies $[x,\mathfrak{s}]\subset \mathfrak{s}_{k}$, then 
$x=0.$ This fact implies the following

\begin{proposition}
The Klein geometry $S/P$ is almost effective with descending filtration
\end{proposition}

\begin{equation}
\{0\}\subset \mathfrak{s}_{k}\subset \mathfrak{p}\subset \mathfrak{s}
\end{equation}

We rewrite (25) as $\mathfrak{s}=\mathfrak{p}^{-}\oplus \mathfrak{s}%
_{0}\oplus \mathfrak{p}^{+}$. Since $\mathfrak{p}^{-},$ $\mathfrak{p}^{+}$
are dual with respect to Killing form, they determine each other and (26)
involves a redundancy. Our purpose is to increase the lentgh of (26) by
exploiting this redundancy.

To simplify things and also to be specific, we assume $S=SL(n,\mathbb{F}),$ $%
n\geq 2,$ $\mathfrak{s}_{0}$ is the Cartan subalgebra of diagonal matrices,
and $\mathfrak{p}$ is the Borel subalgebra of upper triangular matrices so
that $k=n-1 $ in (25). However, the main idea of our construction works in
the broader context of (25).

We start by choosing an abelian subalgebra of $\mathfrak{s}_{-1}$ which we
denote by $\mathfrak{a}_{-1}.$ For instance, some $1$-dimensional subspace
of $\mathfrak{s}_{-1}$ will do. If $n=6,$ then the matrices $A(2,1),$ $%
A(4,3),$ $A(6,5)$ in the standard basis having $1$'s in the indicated
entries and $0$'s elsewhere belong to $\mathfrak{s}_{-1}$ and they commute.
Thus we can choose $\dim \mathfrak{a}_{-1}=3$ in this case. It is easy to
check that $\dim \mathfrak{a}_{-1}$ can be at most $[\frac{n}{2}].$

Having fixed $\mathfrak{a}_{-1},$ we now define $\mathfrak{h}\doteq 
\mathfrak{a}_{-1}\oplus \mathfrak{p}$. Since any $x\in \mathfrak{a}_{-1}$ is
a common eigenvector for all $y\in \mathfrak{s}_{0}$, it follows that $[%
\mathfrak{a}_{-1},\mathfrak{s}_{0}]=\mathfrak{a}_{-1}\subset \mathfrak{h.}$
Using the grading and the fact that $\mathfrak{a}_{-1}$ and $\mathfrak{p}$\
are both subalgebras, we conclude that $\mathfrak{h}\subset \mathfrak{s}$ is
a subalgebra. Note that $\mathfrak{h}$ is neither semisimple nor solvable
(compare to pg. xvi of the Introduction of [11]). We now choose a Lie
subgroup $H\subset S$ with Lie algebra $\mathfrak{h}$ and $H/P$ defines the
descending filtration

\begin{eqnarray}
\mathfrak{h} &=&\mathfrak{a}_{-1}\oplus \mathfrak{p}=\mathfrak{a}_{-1}\oplus 
\mathfrak{s}_{0}\oplus ....\oplus \mathfrak{s}_{n-1} \\
\mathfrak{h}_{0} &=&\mathfrak{p}=\mathfrak{s}_{0}\oplus ....\oplus \mathfrak{%
s}_{n-1}  \notag \\
\mathfrak{h}_{1} &=&\mathfrak{s}_{1}\oplus ....\oplus \mathfrak{s}_{n-1} 
\notag \\
\mathfrak{h}_{2} &=&\mathfrak{s}_{2}\oplus ....\oplus \mathfrak{s}_{n-1} 
\notag \\
&&........  \notag \\
\mathfrak{h}_{n-1} &=&\mathfrak{s}_{n-1}  \notag \\
\mathfrak{h}_{n} &=&0  \notag
\end{eqnarray}

Lemma 4 and Proposition 5 show that (27) is consistent with the order of
jets and (27) shows that $H/P$ is almost effective. We define $G\doteq H/K$
and $G_{0}\doteq P/K$ where $K$ is discrete and obtain

\begin{proposition}
The order of the effective Klein geometry $G/G_{0}$ is $n.$
\end{proposition}

Note that $\mathfrak{h}_{1}/\mathfrak{h}_{2}=\mathfrak{s}_{1}=$ the span of
positive simple roots. Now (22) gives 
\begin{equation}
\begin{array}{ccccccccc}
0 & \longrightarrow & \mathfrak{h}_{1}/\mathfrak{h}_{2} & \longrightarrow & 
\mathfrak{h}/\mathfrak{h}_{2} & \overset{\overline{\pi }_{2,1}}{%
\longrightarrow } & \mathfrak{h}/\mathfrak{h}_{1} & \longrightarrow & 0 \\ 
&  & \downarrow _{\theta _{2}} &  & \downarrow _{\theta _{2}} &  & 
\downarrow _{\theta _{1}} &  &  \\ 
0 & \longrightarrow & J_{2,1}(\overline{\mathfrak{g}})_{o} & \longrightarrow
& J_{2}(\overline{\mathfrak{g}})_{o} & \overset{_{\pi _{2,1}}}{%
\longrightarrow } & J_{1}(\overline{\mathfrak{g}})_{o} & \longrightarrow & 0%
\end{array}%
\end{equation}

Thus $J_{2,1}(\overline{\mathfrak{g}})_{o}$ completely determines the Lie
algebra structure of $\mathfrak{s}.$ However, this is not so for $\mathfrak{h%
}$ due to the arbitrariness involved in the choice of $\mathfrak{a}_{-1}.$

For $\mathbb{F=C}$, $S/P$ is known to be a projective variety. Since $H/P$
is closed in $S/P$ and $G/G_{0}=H/P,$ we obtain

\begin{corollary}
For every integer $n\geq 2$ and $d,$ $1\leq d\leq \lbrack \frac{n}{2}]$,
there exists an effective Klein geometry $G/G_{0}$ of order $n$ such that $%
G/G_{0}$ is a projective variety.
\end{corollary}

We know nothing about the structure of the projective varieties given by
Corollary 9 in general. Even the simplest case $\dim \mathfrak{a}_{-1}=1$
raises some questions about these Riemann surfaces whose answers we do not
know.

We conclude this section with the following remark: The tangent space $%
T(S/P)_{o\text{ }}$of $S/P$ is $\mathfrak{s}/\mathfrak{p}$ $=$ $\mathfrak{s}%
_{-n+1}\oplus \mathfrak{s}_{-n+2}\oplus ....\mathfrak{s}_{-1}$. Therefore
the filtration (27) (starting with $\mathfrak{h}_{1}$) determines a
filtration in $T(S/P)_{o\text{ }}$which, to our knowledge, is observed and
studied first in [2] (see also [3]) in the framework of general parabolic
geometries defined using principal $P$-bundles.

\section{Geometric order}

\bigskip In this section we resume our framework of Section 3, assume that $%
G/G_{0}$ is effective and $G$ is connected. Following [16], we inductively
define $G_{k+1}\doteq \{g\in G_{k}\mid Ad(g)x-x\in \widehat{\mathfrak{g}}_{k}
$ for all $x\in \mathfrak{g},$ $0\leq k\}$ where $\widehat{\mathfrak{g}}_{k}$
is the Lie algebra of $G_{k}.$ Now $G_{k+1}\vartriangleleft G_{k}$ is a
normal subgroup for $k\geq 0$. Thus we obtain the filtration

\begin{equation}
....\subset G_{2}\subset G_{1}\subset G_{0}\subset G
\end{equation}

(29) defines the filtration

\begin{equation}
\mathfrak{....}\subset \widehat{\mathfrak{g}}_{2}\subset \widehat{\mathfrak{g%
}}_{1}\subset \mathfrak{g}_{0}\subset \mathfrak{g}
\end{equation}

where $\widehat{\mathfrak{g}}_{k}$ is the Lie algebra of $G_{k}$ for $k\geq
1.$ Note that $G_{k+1}$ is the kernel of $\ Ad_{k}:G_{0}\rightarrow GL(%
\mathfrak{g}/\widehat{\mathfrak{g}}_{k})$ defined by $Ad_{k}(x)(y+\widehat{%
\mathfrak{g}}_{k})\doteq Ad_{x}(y)+\widehat{\mathfrak{g}}_{k}$ which gives
the faithful representation 
\begin{equation}
Ad_{k}:G_{0}/G_{k+1}\longrightarrow GL(\mathfrak{g}/\widehat{\mathfrak{g}}%
_{k})
\end{equation}

If $\rho :G\rightarrow R$ is \textit{any }homomorphism of Lie groups with
differential $d\rho :\mathfrak{g}\rightarrow \mathfrak{r}$, then $\ker
(d\rho )$ is clearly the Lie algebra of $\ker (\rho ).$ This fact together
with the definitions of $\widehat{\mathfrak{g}}_{k}$ and $\mathfrak{g}_{k}$
shows $\widehat{\mathfrak{g}}_{k}=\mathfrak{g}_{k}.$ In particular, $%
\widehat{\mathfrak{g}}_{m}=\{0\}$ and $G_{m}$ is discrete. If $G/G_{0}$ is
almost effective with $\mathfrak{g}_{m}=\{0\},$ then $K\subset G_{m}$
because (29) stabilizes at $K$ which is therefore discrete$.$

\begin{lemma}
For an effective Klein geometry $G/G_{0}$ with the descending filtration
(9), we have either $i)$ $G_{m}=\{e\}$ or $ii)$ $G_{m}\neq \{e\}$ and $%
G_{m+1}=\{e\}$
\end{lemma}

Proof: Suppose $G_{m}\neq \{e\}.$ Since $G_{m+1}$ is the kernel of $Ad_{m}:$ 
$G_{0}\rightarrow GL(\mathfrak{g}/\widehat{\mathfrak{g}}_{m})=GL(\mathfrak{g}%
),$ we have $G_{m+1}\subset \ker (Ad:G\rightarrow GL(\mathfrak{g}))=Z(G)=$
center of $G$ since $G$ is connected. Therefore $G_{m+1}\subset Z(G)\cap
G_{0}=\{e\}$ since the action is effective. $\ \square $

\begin{definition}
The integer $M$ satisfying $G_{M}=\{e\}$, $G_{M-1}\neq \{e\}$ is called the
geometric order of the effective Klein geometry $G/G_{0}$
\end{definition}

Lemma 10 shows that $M=m$ or $M=m+1$ for an effective Klein geometry. How do
we decide which one is the case? $ii)$ below hints that the answer is
related to the fundamental group of $G/G_{0}.$ This is expected since
fundamental group is the new issue when we pass from Lie algebra to Lie
group.

Setting $k=M-1$ in (31), we obtain

\begin{proposition}
Let $G_{0}$ be any Lie group. Suppose there exists a connected Lie group $G$
satisfying $a)$ $G_{0}\subset G$ is a Lie subgroup $b)$ $G$ acts effectively
on $G/G_{0}.$ Then $G_{0}$ is a matrix group.
\end{proposition}

Now $G_{0}$ is the stabilizer of the point $o$. Therefore, $g\in G_{0}$
defines the isomorphism $\natural j_{k+1}(L_{g})_{o}^{o}:(J_{k}\mathfrak{X)}%
_{o}\rightarrow (J_{k}\mathfrak{X)}_{o}$ and we obtain the representation 
\begin{equation}
Ad_{k}:G_{0}\longrightarrow GL((J_{k}\mathfrak{X)}_{o})
\end{equation}
defined by $g\mapsto \natural j_{k+1}(L_{g})_{o}^{o}$. In fact, for $g\in
G_{0}$, $X\in \mathfrak{g}$, it follows from (7) that 
\begin{eqnarray}
\natural j_{k+1}(L_{g})_{o}^{o}((j_{k}\overline{X})_{o}) &=&j_{k}(\frac{d}{dt%
}(L_{g}\circ L_{e^{tX}}\circ L_{g^{-1}})_{t=0}))_{o}  \notag \\
&=&j_{k}(\frac{d}{dt}(L_{ge^{tX}g^{-1}})_{t=0})_{o}  \notag \\
&=&j_{k}(\overline{Ad(g)X})_{o}
\end{eqnarray}
Note the fundamental difference between (7) and (33): (7) is local whereas
(33) is global and involves analyticity.

(33) shows that (32) restricts as

\begin{equation}
Ad_{k}:G_{0}\longrightarrow GL(J_{k}(\overline{\mathfrak{g}})_{o})
\end{equation}

\begin{lemma}
The kernel of $Ad_{k}$ is $G_{k+1}$.
\end{lemma}

Proof: Let $g\in \ker (Ad_{k})$ so that $\natural j_{k+1}(L_{g})_{o}^{o}$ is
identity on $J_{k}(\overline{\mathfrak{g}})_{o}$, that is, $\natural
j_{k+1}(L_{g})_{o}^{o}((j_{k}\overline{X})_{o})=(j_{k}\overline{X})_{o}$ for
all $X\in \mathfrak{g}_{-1}$. Now (33) gives $j_{k}(\overline{Ad(g)X}%
))_{o}=(j_{k}\overline{X})_{o}$ or equivalently $j_{k}(\overline{Ad(g)X-X}%
)_{o}=0$. By Lemma 4, this condition is equivalent to $Ad(g)X-X\in \mathfrak{%
g}_{k}=\widehat{\mathfrak{g}}_{k}$ and therefore to $g\in G_{k+1}$ by the
definition of $G_{k+1}$. \ $\square $

\begin{proposition}
There is a canonical injection of Lie groups $\Phi
_{k}:G_{0}/G_{k}\rightarrow (\mathfrak{G}_{k})_{o}^{o}\simeq G_{k}(n),$ $%
n=\dim G_{-1}/G_{0}$.
\end{proposition}

Proof: We define the homomorphism $G_{0}\rightarrow (\mathfrak{G}%
_{k})_{o}^{o}$ by $g\rightarrow j_{k}(L_{g})_{o}^{o}.$ By Lemma 1, $%
j_{k}(L_{g})_{o}^{o}=id$ if and only if $\natural j_{k}(L_{g})_{o}^{o}=id$
and the conclusion follows from Lemma 13. \ $\square $

Lemma 13 gives the faithful representation $Ad_{k}:G_{0}/G_{k+1}\rightarrow
GL(J_{k}(\overline{\mathfrak{g}})_{o}),$ and Propositions 5, 14 give the
commutative diagram

\begin{equation}
\begin{array}{cccc}
Ad_{k}: & G_{0}/G_{k+1} & \longrightarrow & GL(\mathfrak{g}/\mathfrak{g}_{k})
\\ 
& \parallel \Phi _{k+1} &  & \parallel \theta _{k} \\ 
\natural _{\mid (G,G_{0})}: & \Phi _{k+1}(G_{0}/G_{k+1}) & \longrightarrow & 
GL(J_{k}(\overline{\mathfrak{g}})_{o})%
\end{array}%
\end{equation}
which is the group analog of (24).

Setting $k=M$ in Proposition 14 gives

\begin{corollary}
There is a canonical injection of Lie groups $\Phi _{M}:G_{0}\rightarrow $ $(%
\mathfrak{G}_{M})_{o}^{o}$
\end{corollary}

The following two special cases of Corollary 15 are of special interest.

$i)$ Let $P\subset SL(n,\mathbb{C})$ be the Borel subgroup of upper
triangular matrices. Then there exists a discrete and normal subgroup $%
K\subset P$ such that $P/K$ injects canonically into $(\mathfrak{G}%
_{1})_{o}^{o}\simeq G_{M}(1),$ where $M=n$ or $n+1.$ To prove this
statement, we recall the construction of $G/G_{0}$ in Proposition 8 and
choose $\dim \mathfrak{a}_{-1}=1.$

Before we state $ii),$ we recall the main result of [15]: A smooth manifold $%
X$ is determined up to diffeomorphism by the abstract Lie algebra structure
of $\mathfrak{X.}$ This statement holds also in the analytic category and we
refer to [5] for an extensive literature on the generalizations of this
result. It turns out that certain subalgebras of $\mathfrak{X}$ also
determine $X.$ In the analytic category, all the information in $\mathfrak{X}
$ is encoded at one point. The next statement shows how the fundamental
group is encoded in first order jets in a special case. \ 

$ii)$ Let $G$ be simply connected, $G_{0}\subset G$ discrete and $G/G_{0}$
effective. Then we have the canonical injection $\pi
_{1}(G/G_{0},o)\rightarrow (\mathfrak{G}_{1})_{o}^{o}$ and the faithful
representation $\pi _{1}(G/G_{0},o)\rightarrow GL(\mathfrak{g}).$ For the
proof, we note that $G_{0}\simeq \pi _{1}(G/G_{0},o)$ since $G_{0}$ is
simply connected and $m=0$ since $G_{0}$ is discrete. Assuming that $%
G_{0}\neq \{e\},$ we conclude $M\neq 0$ and therefore $M=1$ by Lemma 10. We
now set $k=M-1$ in (35).

The next corollary is an extension of Corollary 15.

\begin{corollary}
For any integer $k\geq 1,$ there is a canonical injection $\Phi
_{M+k}:G_{0}\rightarrow (\mathfrak{G}_{M+k})_{o}^{o}$ which makes the
following diagram commute:
\end{corollary}

\begin{equation}
\begin{array}{ccc}
(\mathfrak{G}_{M+k})_{o}^{o} & \overset{\pi _{M+k,M}}{\longrightarrow } & (%
\mathfrak{G}_{M})_{o}^{o} \\ 
& \nwarrow _{\Phi _{M+k}} & \uparrow _{\Phi _{M}} \\ 
&  & G_{0}%
\end{array}%
\end{equation}

Proof: We define $\Phi _{M+k}$ as in the proof of Proposition 14. $\square $

In the same way as Levi-Civita connection is the object canonically
determined by a Riemannian structure, the splitting given by (36) implies
the existence of some objects canonically determined by some geometric
structures. This problem will be studied elsewhere (see [2], [3] for the
standard but very different approach to this problem).

We single out the next fact, which is essentially equivalent to Corollary
15, as a corollary for reference in Section 6.

\begin{corollary}
Let $a,b\in G$ satisfy $L_{a}(p)=L_{b}(p)=q$, for some $p,q\in G/G_{0}$. If $%
j_{M}(L_{a})_{q}^{p}=j_{M}(L_{b})_{q}^{p}$, then $a=b$.
\end{corollary}

Proof: First we assume $p=o$. Now $j_{M}(L_{a})_{q}^{o}=j_{M}(L_{b})_{q}^{o}$
iff $j_{M}(L_{a^{-1}})_{o}^{q}\circ
j_{M}(L_{b})_{q}^{o}=j_{M}(L_{a^{-1}b})_{o}^{o}=$ $id.$ Therefore $%
a^{-1}b\in G_{M}=\{e\}$ and $a=b$. The claim for arbitrary $p$ follows from
homogeneity. \ $\square $

\section{Differential equations of finite type}

In this section, we derive an important consequence of Corollary 17.

Following Lie, we write the action of $G$ on $G/G_{0}$ locally as

\begin{equation}
f^{i}(x^{1},x^{2},...,x^{r},y^{1},y^{2},...,y^{n})=z^{i}\quad 1\leq i\leq n
\end{equation}
where $\dim G/G_{0}=n$ and $(y^{i}),$ $(z^{i})$ are local parameters for $%
G/G_{0},$ $\dim G=r$ and $(x^{i})$ are local parameters for $G.$ We write
(37) shortly as $f(x,y)=z$ and fix some $\overline{x},\overline{y},\overline{%
z}$ with $f(\overline{x},\overline{y})=\overline{z}$. Thus $\overline{x}\in
G $ determines the diffeomorphism $z=z(y)$ defined by $f(\overline{x},y)=z$.
The $k$-jet $j_{k}(\overline{x})_{\overline{z}}^{\overline{y}}$ of this
diffeomorphism is given by

\begin{equation}
\frac{\partial ^{\left| \mu \right| }f^{i}(\overline{x},\overline{y})}{%
\partial y^{\mu }}=\frac{\partial ^{\left| \mu \right| }z^{i}}{\partial
y^{\mu }}(\overline{y})\qquad 0\leq \left| u\right| \leq k
\end{equation}

Corollary 17 asserts that $\overline{x}$ is uniquely determined by (38) for $%
0\leq \left| u\right| \leq M.$ Thus we can solve $\overline{x}$ in terms of $%
\overline{y}$ and $\frac{\partial ^{\left| \mu \right| }z}{\partial y^{\mu }}%
(\overline{y})$ from (38) for $0\leq \left| u\right| \leq M$ and substitute
the result into (38) for $0\leq \left| u\right| \leq M+1.$ The result is 
\begin{equation}
\Phi _{\mu }^{i}(\overline{y},\frac{\partial ^{\left| \sigma \right| }z^{i}}{%
\partial y^{\sigma }}(\overline{y}))=\frac{\partial ^{\left| u\right| }z^{i}%
}{\partial y^{\mu }}(\overline{y})\quad 1\leq \left| u\right| \leq M+1,\text{
}0\leq \left| \sigma \right| \leq M
\end{equation}
for some functions $\Phi _{\mu }^{i}.$ The form of (39) does not depend on
the choices $\overline{x},\overline{y}$. Regarding $\overline{y}$ as a
variable, (39) defines a system of $PDE$'s of finite type. If $\dim
G/G_{0}=1,$ then (39) is a system of $ODE$'s. We can write (39) in a
coordinate system such that we have also $0\leq \left| \sigma \right| \leq
\left| u\right| .$ The number of equations in (39) is $\dim G_{M+1}(n)$ and $%
\dim G_{0}$ of them are dependent.

To summarize, we have

\begin{proposition}
An effective Klein geometry $G/H$ of geometric order $M$ defines a
differential equation on $G/H$ of order $M+1.$ The global solution space of
this differential equation is $G.$
\end{proposition}

In order to solve the group parameters uniquely from (38), the number of
unknowns should not exceed the number of equations. This gives an inequality
relating $M$, $\dim G$ and $\dim H,$ which is analogous to the inequalities
in [12] (see pg.161-162).

Two instances of Proposition 18 are well known:

$1)$ Note that $M=0$ iff $G_{0}=\{e\}$ so that a Klein geometry $G/G_{0}$ of
geometric order zero is nothing but the Lie group $G$ together with the left
action of $G$ on itself. In view of the simple computation in the proof of
Theorem 5 in [14] on page 178, (39) becomes

\begin{equation}
\frac{\partial z^{i}}{\partial y^{j}}=\xi _{a}^{i}(z)\omega _{j}^{a}(y)
\end{equation}
where the vector fields $\xi _{j}=\xi _{j}^{a}\partial _{a}$ are the
infinitesimal generators of the action and the 1-forms $\omega ^{i}=\omega
_{a}^{i}dx^{a}$ on $G$ are components of the Maurer-Cartan form. Therefore
Proposition 18 reduces to Lie's First Fundamental Theorem.

$2)$ Let $G=SL(2,\mathbb{C)}$ and $H=$ the subgroup of upper trangular
matrices and $K=\{\pm I\}.$ Now $G\doteq G/K$ can identified with the group $%
\mathcal{M}$ of normalized M\"{o}bius transformations

\begin{equation}
\frac{az+b}{cz+d}=w,\quad ad-bc=1,
\end{equation}
$G/H=\frac{G/K}{H/K}\doteq G/G_{0}$ with the complex sphere $\mathcal{S}$
and the effective action of $G$ on $G/G_{0}$ with the effective action of $%
\mathcal{M}$ on $\mathcal{S}$ as point transformations. In this case $M=2$
and the process of deriving (39) from (37) amounts to differentiating (41)
three times and eliminiting the group parameters (see pg.21 of [14] for
details). We have $\dim G_{3}(1)=3$ equations in (39) and $\dim G_{0}=2$ of
these equations are dependent. In fact, the first two equations reduce to
identities and the third one is the well known Schwarzian $ODE$

\begin{equation}
w^{\prime \prime \prime }=\frac{3}{2}\frac{(w^{\prime \prime })^{2}}{%
w^{\prime }}
\end{equation}

We refer to [14], [12] on the relation of (42) to differential invariants
and to [9] for the use of (42) in defining projective structures on Riemann
surfaces.

Note that the second statement of Proposition 18 is rather a definition than
an assertion. We refer to [19] for a categorical approach to the global
formulation of $PDE$'s and their symmetries in terms of diffieties and to
the classical book [11] for a systematic study of symmetries of differential
equations.

We arrived at (39) starting from geometry. Conversely, we can also start
with differential equations of finite type and study their geometrization as
in [21], [22]. It is therefore no coincidence that semisimple Lie groups,
parabolic subgroups and projective imbeddings arise naturally also in [21],
[22] as in this paper. This geometrization problem is studied for general
exterior differential systems in the influential paper [18].

\bigskip It is standard to take the bundle $G\rightarrow G/H$ as the basis
of a Klein geometry and generalize $G$ to an auxiliary but extremely useful
object on which a Lie group $H$ acts freely (see the Foreword of [16] by
S.S.Chern for a very concise formulation of this point of view). It is well
known that this approach has been remarkably successful with far reaching
results. However, we believe that the approach to geometry based on
transitive actions of Lie groups has also much to offer.

\section{Appendix}

\bigskip All ingredients of our constructions in Section 3 are contained in
the fundamental papers [8], [17] in much more generality and we make here
some comments on the relation of Section 3 to these works, emphasizing the
novelty of our approach.

Let $\mathfrak{g}\subset gl(V)$ be a subalgebra. Following [8], [17], [7],
we set $pr_{0}\mathfrak{g}\doteq \mathfrak{g}$ and define the $k$'th
prolongation $pr_{k}\mathfrak{g}$ of $\mathfrak{g}$ inductively by the
formula $pr_{k}\mathfrak{g}\doteq \{S\in Hom(V,pr_{k-1}\mathfrak{g}$ $)\mid
S(v)w=S(w)v,$ for all $v,w\in V\}.$

\begin{definition}
$\mathfrak{g}\subset gl(V)$ is called finite type if $pr_{m}\mathfrak{g}=0$
for some $m.$ If $pr_{\widetilde{m}}\mathfrak{g}=0$ and $pr_{\widetilde{m}-1}%
\mathfrak{g}\neq 0$, $\widetilde{m}$ is called the prolongation order of $%
\mathfrak{g}\subset gl(V)$.
\end{definition}

The best known matrix algebras which are finite type are the orthogonal,
conformal, projective and affine algebras (see [8], [17] for details). As
far as we know, for all first order $G$-structures studied so far in
geometry, either the Lie algebra $\mathfrak{g}$ of $G$ is of infinite type
(as in symplectic or complex structures) or of finite type with $\widetilde{m%
}\leq 2.$

Next setting $k=0$ in (24), we formulate

\begin{definition}
$ad_{0}(\mathfrak{g}_{0}/\mathfrak{g}_{1})\subset gl(\mathfrak{g}/\mathfrak{g%
}_{0})$ is called the first order isotropy algebra of the effective Klein
geometry $G/G_{0}.$
\end{definition}

It is now easy to prove the following

\begin{proposition}
The infinitesimal order $m$ of the first order isotropy algebra $ad_{0}(%
\mathfrak{g}_{0}/\mathfrak{g}_{1})\subset gl(\mathfrak{g}/\mathfrak{g}_{0})$
of an effective Klein geometry $G/G_{0}$ satisfies $m\leq \widetilde{m}$
(where we set $\widetilde{m}=\infty $ if $ad_{0}(\mathfrak{g}_{0}/\mathfrak{g%
}_{1})$ is of infinite type).
\end{proposition}

Also, the concept of a rigid geometric structure is introduced in [6]. These
structures are characterized by the fact that their infinitesimal
automorphisms are determined by their jets of some fixed order. In view of
our results in this paper, it is not surprising that this concept turns out
to be equivalent to first order $G$-structures being of finite type, as
shown in [1]. We believe that there is much conceptual overlap with [6], [1]
and this paper.\bigskip 

\bigskip

\bigskip

\textbf{Acknowledgements:} Sections 5, 6 did not exist in the original
version of this paper. We are indebted to P.J.Olver for drawing our
attention to [21], [22] which gave birth to Proposition 18. Part of this
work was completed when the second author was visiting Feza G\"{u}rsey
Institute. He is grateful to T\"{u}bitak for its financial support and to
T.Turgut for his encouragement for this research.

\bigskip

\bigskip

\bigskip

\bigskip \textbf{References}

\bigskip

[1] A.Candel, R.Quiroga-Barranco: Gromov's centralizer theorem, Geom.
Dedicata 100 (2003), 123-135

[2] A.Cap, H.Schichl: Parabolic geometries and canonical Cartan connections,
Hokkaido Math. J. 29. No.3 (2000), 453-505

[3] A.Cap, J.Slovak, V.Soucek: Bernstein-Gelfand-Gelfand Sequences, Ann. of
Math., 154 (2001), 97-113

[4] D.Fuks: Cohomology of Infinite Dimensional Lie Algebras, Contemporary
Soviet Mathematics, Consultants Bureau, New York and London, 1986

[5] J.Grabowski, N.Poncin: Lie algebraic characterization of manifolds,
Cent. Eur. J. Math., 2(5), 2005, 811-825

[6] M.Gromov: Rigid transformation groups, Geometrie differentielle (Paris,
1986), 65-139, Travaux en Cours, 33, Hermann, Paris, 1988

[7] V.W.Guillemin: The integrability problem for G-structures, Trans.Amer.
Math.Soc. 116 1965 544-560

[8] V.W.Guillemin, S.Sternberg: An algebraic model of transitive
differential geometry, Bull. Amer. Math. Soc. 70 1964, 16-47

[9] R.C.Gunning: Lectures on Riemann Surfaces, Princeton Mathematical Notes,
Princeton University Press, Princeton, N.J. 1966

[10] I.S.Krasilschik, V.V.Lychagin, A.M.Vinogradov: Geometry of Jet Spaces
and Nonlinear Partial Differential Equations, Translated from the Russian by
A. B. Sosinski. Advanced Studies in Contemporary Mathematics, 1. Gordon and
Breach Science Publishers, New York, 1986

[11] P.J.Olver: Applications of Lie groups to Differential Equations,
Graduate Texts in Mathematics, Springer-Verlag, New York Inc., 1986

[12] P.J.Olver: Equivalence, Invariants, and Symmetry, Cambridge University
Press, 1995

[13] E.Orta\c{c}gil: The heritage of S.Lie and F.Klein: Geometry via
transformation groups, arXiv.org, math.DG/0604223, posted on April 2006

[14] J.F.Pommaret: Partial Differential Equations and Group Theory, New
Perspectives for Applications, Kluwer Academic Publishers, Vol. 293, 1994

[15] M.E.Shanks, L.E.Pursell: The Lie algebra of a smooth manifold, Proc,
Amer. Math. Soc. 5, (1954), 468-472

[16] R.W.Sharpe: Differential geometry. Cartan's generalization of Klein's
Erlangen program, Graduate Texts in Mathematics, 166. Springer-Verlag, New
York, 1997

[17] I.M.Singer, S.Sternberg: The infinite groups of Lie and Cartan. I. The
transitive groups. J. Analyse Math. 15 1965 1--114

[18] N.Tanaka: On the equivalence problems associated with simple graded Lie
algebras, Hokkaido Math. J. 8 (1979), 23-84

[19] A.M.Vinogradov: Cohomological Analysis of Partial Differential
Equations and Secondary Calculus, translated from the Russian manuscript by
Joseph Krasilshik, Translations of Mathematical Monographs, 204, American
Mathematical Society, Providence, RI, 2001 (Reviewer: Victor V. Zharinov)

[20] G.Weingart: Holonomic and semi-holonomic geometries, Global Analysis
and harmonic analysis (Marseille-Luminy, 1999), 307-328, Semin. Congr., 4,
Soc. Math. France, Paris, 2000

[21] K.Yamaguchi, T.Yatsui: Geometry of higher order differential equations
of finite type associated with symmetric spaces, Advanced Studies in Pure
mathematics, Lie Groups, Geometric Structures and Differential Equations,
-One Hundred Years after Sophus Lie-, 0-61, 2000

[22] K.Yamaguchi, T.Yatsui: Parabolic geometries associated with
differential equations of finite type, preprint

\bigskip

Ender Abado\u{g}lu, Yeditepe \"{U}niversitesi, Matematik B\"{o}l\"{u}m\"{u},
26 A\u{g}ustos Yerle\c{s}imi, 81120, Kay\i \c{s}da\u{g}\i , \.{I}stanbul, T%
\"{u}rkiye, e-mail: eabadoglu@yeditepe.edu.tr

\bigskip

Erc\"{u}ment Orta\c{c}gil, Bo\u{g}azi\c{c}i Universitesi, Matematik B\"{o}l%
\"{u}m\"{u}, 34342, Bebek, \.{I}stanbul, T\"{u}rkiye, e-mail:
ortacgil@boun.edu.tr

\bigskip

Ferit \"{O}zt\"{u}rk, Bo\u{g}azi\c{c}i Universitesi, Matematik B\"{o}l\"{u}m%
\"{u}, 34342, Bebek, \.{I}stanbul, T\"{u}rkiye, e-mail:
ferit.ozturk@boun.edu.tr

\end{document}